\documentclass[11pt]{article}
\usepackage[margin=1in]{geometry}
\usepackage[authoryear,round]{natbib}
\usepackage{amsmath,amssymb,amsthm}
\usepackage{graphicx}
\usepackage{booktabs}
\usepackage{hyperref}
\usepackage{thm-restate}

\newtheorem{lemma}{Lemma}
\theoremstyle{remark}
\newtheorem{remark}{Remark}

\newcommand{\doi}[1]{\href{https://doi.org/#1}{doi:\nolinkurl{#1}}}
\newcommand{\nuB}{\nu_{\mathrm{B}}}
\newcommand{\nuG}{\nu_{\Gamma}}
\DeclareMathOperator{\logit}{logit}

\title{A simple improvement on Kerman's approximation\\to the median of the beta distribution}
\author{Richard F. Lyon\\[4pt]
{\normalsize Google Research, Mountain View, California, USA}\\
{\normalsize \texttt{dicklyon@acm.org}}}
\date{\today}

\begin{document}
\maketitle

\begin{abstract}
Kerman's approximation to the median of the $\mathrm{Beta}(\alpha,\beta)$ distribution is the
ratio $g(\alpha)/(g(\alpha) + g(\beta))$ with $g(x) = x - \frac13$, using the first two terms of the
asymptotic expansion of the median of the standard gamma distribution. We extend this elegant ratio
form to the one-parameter family with $g(x) = x - \frac13 + c_1/x$, in which $c_1 = 0$ gives Kerman's.
For large shape parameters along any ray $\beta = r\alpha$ with $r > 1$, the resulting
approximation is eventually a lower bound on the median if $c_1 \le \frac{8}{405}$ and an upper
bound if $c_1 \ge \frac{1}{45}$. We recommend $c_1 = \frac{1}{36}$ for its simplicity: it makes $g(x) = (x - \frac16)^2/x$
a perfect square over $x$, and it is within $1\%$ of the minimax value of $c_1$ for $\alpha, \beta \ge 1$
under several error measures, while the theoretically distinguished values $\frac{8}{405}$ and $\frac{1}{45}$
are less accurate. The approximation is then
\[
\nu_{1/36}(\alpha,\beta) = \frac{(\alpha - \frac16)^2/\alpha}{(\alpha - \frac16)^2/\alpha + (\beta - \frac16)^2/\beta} \, .
\]
Its relative error lies between $-0.171\%$ and $+0.212\%$ for all $1 \le \alpha \le \beta$, compared
with up to $-3.8\%$ for Kerman's. It is more accurate than Kerman's approximation wherever the
latter is defined, with relative error never more than $\frac{13}{32}$ of Kerman's, the ratio approached
asymptotically. The same $g(\alpha)$ is a useful approximation to the gamma distribution median,
overestimating it by at most $0.212\%$ for shape parameter $\alpha \ge 1$.

\end{abstract}

\section{Introduction}

The median $\nuB(\alpha,\beta)$ of the beta distribution with shape parameters $\alpha$ and $\beta$, the
solution $x$ of $I_x(\alpha,\beta) = \frac12$, has no closed form except in special cases. Simple
approximations to it are useful because it comes up often: it is the posterior median of a binomial
proportion under a conjugate beta prior, such as the neutral prior $\mathrm{Beta}(\frac13, \frac13)$ of
\citet{kerman2011neutral}, whose posteriors have both shape parameters at least $\frac13$; it is the median of the $k$th order statistic of a uniform sample
of size $n$, which has the $\mathrm{Beta}(k, n-k+1)$ distribution, and so gives the median ranks used in
probability plotting; and it is a natural starting value for the numerical inversion of $I_x$.

\citet{kerman2011closed} proposed the approximation
\begin{equation}
\nuB(\alpha,\beta) \approx \frac{\alpha - \frac13}{\alpha + \beta - \frac23} \, ,
\label{eq:kerman}
\end{equation}
which is asymptotically accurate as both shape parameters grow. It has the \emph{ratio form}
$g(\alpha)/(g(\alpha) + g(\beta))$ with $g(x) = x - \frac13$. We study the one-parameter family
\begin{equation}
g_{c_1}(x) = x - \frac13 + \frac{c_1}{x} \, , \qquad
\nu_{c_1}(\alpha,\beta) = \frac{g_{c_1}(\alpha)}{g_{c_1}(\alpha) + g_{c_1}(\beta)} \, ,
\label{eq:family}
\end{equation}
with $c_1 \ge 0$, in which $c_1 = 0$ gives Kerman's approximation. Our recommended member has
$c_1 = \frac{1}{36}$, for which $g$ factors as $g_{1/36}(x) = (x - \frac16)^2/x$, and $\nu_{1/36}$ can be written in
several forms:
\begin{equation}
\begin{split}
\nu_{1/36}(\alpha,\beta) &= \frac{g_{1/36}(\alpha)}{g_{1/36}(\alpha) + g_{1/36}(\beta)}
= \frac{(\alpha - \frac16)^2/\alpha}{(\alpha - \frac16)^2/\alpha + (\beta - \frac16)^2/\beta} \\[4pt]
&= \frac{\alpha - \frac13 + \frac{1}{36\alpha}}
{\alpha + \beta - \frac23 + \frac{1}{36\alpha} + \frac{1}{36\beta}} \\[4pt]
&= \frac{(\alpha - \frac16)^2\,\beta}{(\alpha - \frac16)^2\,\beta + (\beta - \frac16)^2\,\alpha}
= \frac{(6\alpha - 1)^2\,\beta}{(6\alpha - 1)^2\,\beta + (6\beta - 1)^2\,\alpha} \, .
\end{split}
\label{eq:nuB3}
\end{equation}
It is nearly as simple as Kerman's approximation but much more accurate: over $1 \le \alpha \le \beta$ its
relative error is at most $0.212\%$, against $3.82\%$ for Kerman's, and $\frac{1}{36} \approx 0.027778$
is within $1\%$ of the minimax value of $c_1$ there, $0.027537$, whose worst error is $0.189\%$.
Since $g_{1/36}$ has a double zero at $x = \frac16$ and grows again below it, $\nu_{1/36}$ is intended for
$\alpha, \beta > \frac16$; it is most accurate for $\alpha, \beta \ge 1$, and its relative error stays within $3.3\%$ for
$\alpha, \beta \ge \frac12$ and $14.4\%$ for $\alpha, \beta \ge \frac13$ (Table~\ref{tab:errors}).

Section~\ref{sec:background} reviews the properties of ratio forms, the gamma median expansion behind
them, and two older approximations of the same kind. Section~\ref{sec:rays} locates $\nu_{1/36}$ within
the family by its asymptotic behavior for large shape parameters. Section~\ref{sec:nuB3} compares its
accuracy with Kerman's and other members numerically. Section~\ref{sec:idealc} maps the value of $c_1$
that would make the ratio form exact at each $(\alpha,\beta)$. Section~\ref{sec:otherc} discusses
minimax and shape-dependent choices of $c_1$. Finally, Section~\ref{sec:gamma} briefly explores
corresponding approximations to the gamma median.

\section{Background and related approximations}
\label{sec:background}

Every ratio form is exact on the diagonal $\alpha = \beta$ of the $(\alpha,\beta)$ plane, where the
median is $\frac12$, and respects the symmetry $\nuB(\beta,\alpha) = 1 - \nuB(\alpha,\beta)$. When $\beta$
is large, the beta distribution approaches a gamma distribution scaled by $1/\beta$, so
$\beta\,\nuB(\alpha,\beta) \to \nuG(\alpha)$, the median of the gamma distribution with shape parameter
$\alpha$ and unit scale; a ratio form with $g(x) \sim x$ correspondingly has
$\beta\,g(\alpha)/(g(\alpha)+g(\beta)) \to g(\alpha)$. So a good gamma median approximation $g$ makes a good
beta median approximation at large $\beta$. The gamma median has no closed form either, but it has a
well-known asymptotic expansion \citep{choi1994medians, berg2006chen},
\begin{equation}
\nuG(\alpha) = \alpha - \frac13 + \frac{8}{405\alpha} + \frac{184}{25515\alpha^2} + \cdots \, ,
\label{eq:gammaseries}
\end{equation}
and $\alpha - \frac13$ is a lower bound on it \citep{chen1986bounds, berg2006chen}. Kerman's $g$ consists of the
first two terms of \eqref{eq:gammaseries}, and truncating after the third term gives $\nu_{8/405}$.

An older ratio form is implicit in the normal approximation to the $F$ distribution of
\citet{paulson1942approximate}, which applies the cube-root transformation of \citet{wilson1931distribution}
to the numerator and denominator of $F$. Setting Paulson's normal deviate to zero gives the median of $F$
with $n_1$ and $n_2$ degrees of freedom as $\big((1 - \frac{2}{9n_1})/(1 - \frac{2}{9n_2})\big)^3$.
Since $\alpha F/(\alpha F + \beta)$ has the $\mathrm{Beta}(\alpha,\beta)$ distribution when $n_1 = 2\alpha$
and $n_2 = 2\beta$, this is the ratio form with
\[
g_{\mathrm{WH}}(x) = x\Big(1 - \frac{1}{9x}\Big)^3 = x - \frac13 + \frac{1}{27x} - \frac{1}{729x^2} \, .
\]
The resulting $\nu_{\mathrm{WH}}$ is already more accurate than Kerman's approximation: over
$1 \le \alpha \le \beta$ it overestimates the median by at most $1.33\%$, where Kerman's underestimates it by
up to $3.82\%$, and it is better on every domain of Table~\ref{tab:errors}. But it is little known as an
approximation to the beta median, since it appears only implicitly, as the median case of an approximation
to the $F$ distribution.
Kerman's $g$, $g_{\mathrm{WH}}$ and the series \eqref{eq:gammaseries} all have the form
$x - \frac13 + c_1/x + c_2/x^2 + \cdots$: Kerman's with $c_1 = c_2 = 0$, $g_{\mathrm{WH}}$ with
$c_1 = \frac{1}{27} = \frac{15}{405}$ and $c_2 = -\frac{1}{729}$, and the series with
$c_1 = \frac{8}{405}$ and $c_2 = \frac{184}{25515}$. Section~\ref{sec:rays} shows that for large shape
parameters only $c_1$ matters to second order, which is why we concentrate on the family \eqref{eq:family}.

Approximations to percentage points of Fisher's $z = \frac12\log F$ also give the beta median, but they
are designed for the tails. At the median, \citeauthor{carter1947approximation}'s
\citeyearpar{carter1947approximation} refinement of the Fisher--Cochran formula reduces to
$\logit\nu = \log(\alpha/\beta) - \frac23(u - v)(1 - u - v)$, with $u = 1/(2\alpha - 1)$ and
$v = 1/(2\beta - 1)$. It is singular at $\alpha$ or $\beta = \frac12$, and to second order its logit matches
that of a ratio form with $c_1 = \frac{1}{18}$; over $1 \le \alpha \le \beta$ its worst relative error is
$+44\%$, at $\alpha = 1$ as $\beta \to \infty$, where it gives $1/\beta$ in place of $(\log 2)/\beta$.

\section{Asymptotic bounds along rays}
\label{sec:rays}

When both shape parameters are large, every member of the family \eqref{eq:family} agrees with the
median to first order in $1/\alpha$; the members differ at second order, where $c_1$ alone decides whether
the approximation lies above or below the median. Along each ray $\beta = r\alpha$ with $r > 1$, the dividing
value of $c_1$ is explicit. The two theorems below are proved in Appendix~\ref{app:proofs}.

\begin{restatable}[Asymptotic bounds]{theorem}{thmbounds}
\label{thm:cfamily}
Fix $r > 1$ and let
\[
c^*_r = \frac{8}{405} + \frac{4r}{405\,(1+r)^2} \, .
\]
Let $g(x) = x - \frac13 + c_1/x + O(x^{-2})$ as $x \to \infty$. If $c_1 < c^*_r$, then
$g(\alpha)/(g(\alpha) + g(r\alpha)) < \nuB(\alpha, r\alpha)$ for all sufficiently large $\alpha$.
If $c_1 > c^*_r$, the reverse inequality holds for all sufficiently large $\alpha$.
\end{restatable}

The threshold has a simple meaning: along each ray, the median itself agrees to second order with the
ratio form $\nu_{c^*_r}$, as \eqref{eq:logiterr} in the proof shows.
Since $0 < 4r/(1+r)^2 < 1$ for $r > 1$, we have $\frac{8}{405} < c^*_r < \frac{9}{405} = \frac{1}{45}$.
So $c_1 \le \frac{8}{405}$ gives an asymptotic lower bound on every ray, $c_1 \ge \frac{1}{45}$ an asymptotic
upper bound on every ray, and $c_1$ strictly between gives one or the other depending on the ray.
In particular, Kerman's $\nu_0$, the truncated series $\nu_{8/405}$ and the ratio form of the true
gamma median itself are asymptotic lower bounds, while $\nu_{1/36}$ is an asymptotic upper bound.

\begin{remark}
\label{rem:threshold}
At $c_1 = \frac{1}{45}$ the second-order term vanishes as $r \to 1$, so near the diagonal the sign is decided
at third order. Numerically, $c_1 = \frac{1}{45}$ fails to be an upper bound in a thin region near the
diagonal ($\beta/\alpha < 1.51$ at $\alpha = 100$) that shrinks as $\alpha$ grows.
\end{remark}

\begin{restatable}[Asymptotic error]{theorem}{thmerror}
\label{prop:ratio}
Fix $r > 1$ and $c_1 \ge 0$. As $\alpha \to \infty$ with $\beta = r\alpha$, the relative error
$e_{c_1} = \nu_{c_1}/\nuB - 1$ satisfies
\[
e_{c_1}(\alpha, r\alpha) = -\Big(1 - \frac1r\Big)\,\frac{c^*_r - c_1}{\alpha^2} + O(\alpha^{-3}) \, .
\]
In particular
\[
\frac{e_{1/36}}{e_0} \to -\frac{\frac{1}{36} - c^*_r}{c^*_r} \, ,
\]
whose magnitude increases from $\frac14$ as $r \to 1$ to $\frac{13}{32}$ as $r \to \infty$.
\end{restatable}

For example, along the ray $\beta = 2\alpha$, Kerman's approximation has relative error about
$-0.0110/\alpha^2$ and $\nu_{1/36}$ about $+0.0029/\alpha^2$.

\section{The recommended approximation with \texorpdfstring{$c_1 = \frac{1}{36}$}{c = 1/36}}
\label{sec:nuB3}

\begin{figure}[t]
  \centering
  \includegraphics{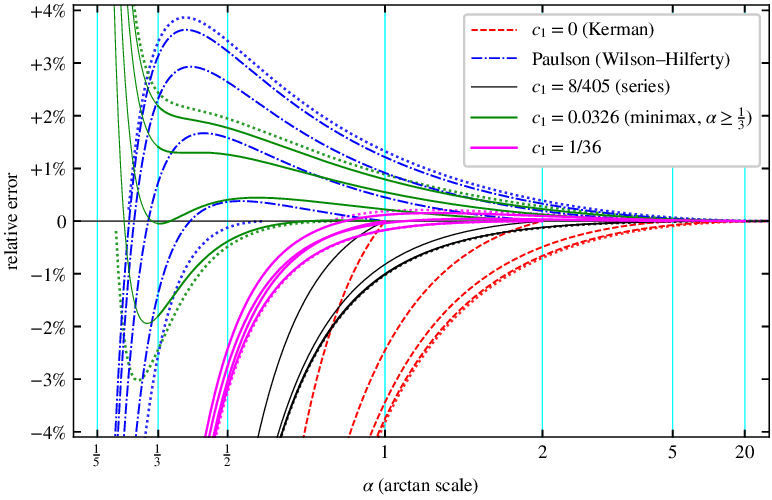}
  \caption{Relative errors of ratio-form approximations to the beta median for $\beta = 1$, $2$, $5$
  and $20$, each curve ending at $\alpha = \beta$ where all are exact, on an $\arctan\alpha$ axis.
  Kerman's $\nu_0$ and the truncated asymptotic series $\nu_{8/405}$ fall below the median;
  $\nu_{1/36}$ stays within $\pm 0.22\%$ for $\alpha \ge 1$, and Paulson's $\nu_{\mathrm{WH}}$ (blue dash-dot)
  stays above it there, by up to $1.33\%$. The green curves use the $c_1$ that
  minimizes the worst error over $\alpha \ge \frac13$; they are drawn thin below $\alpha = \frac13$, outside
  that domain, where the $c_1/x$ term turns them sharply upward as $\alpha$ approaches the minimum of $g$
  at $\sqrt{c_1} \approx 0.18$. Dotted lines of the same colors are the envelopes, the most negative
  and most positive errors over all $\beta \ge \alpha$ including $\beta \to \infty$.}
  \label{fig:rel}
\end{figure}

Figure~\ref{fig:rel} compares the relative errors of $\nu_{1/36}$ and other members of the family for
several values of $\beta$, and Table~\ref{tab:errors} gives the worst errors over
three domains. On $1 \le \alpha \le \beta$, the relative error of $\nu_{1/36}$ lies between
$-0.171\%$ (near $\alpha = 1$, $\beta \approx 1.84$) and $+0.212\%$ (near $\alpha = 1.18$ as
$\beta \to \infty$), an eighteenth of Kerman's worst error there.

\begin{table}[t]
\centering
\caption{Worst negative and positive relative errors (in percent) of ratio forms with
$g(x) = x - \frac13 + c_1/x + c_2/x^2$ over $\alpha_0 \le \alpha \le \beta$, computed on a grid
(Section~\ref{sec:numerics}). A zero means no error of that sign beyond rounding ($10^{-6}\%$) was found.
The second and third rows use the series \eqref{eq:gammaseries} truncated after its $c_1$ and $c_2$ terms,
and the last row uses Paulson's $g_{\mathrm{WH}}$.}
\label{tab:errors}
\medskip
\small
\setlength{\tabcolsep}{5pt}
\renewcommand{\arraystretch}{1.3}
\begin{tabular}{lccccc}
\toprule
Approximation & $c_1$ & $c_2$ & $\alpha_0 = \frac13$ & $\alpha_0 = \frac12$ & $\alpha_0 = 1$ \\
\midrule
Kerman, $\nu_0$ & $0$ & $0$ & $[-100,\ 0]$ & $[-26.73,\ 0]$ & $[-3.820,\ 0]$ \\
series, to $c_1$ & $\frac{8}{405}$ & $0$ & $[-38.05,\ 0]$ & $[-9.448,\ 0]$ & $[-1.020,\ 0]$ \\
series, to $c_2$ & $\frac{8}{405}$ & $\frac{184}{25515}$ & $[-0.191,\ +29.99]$ & $[-0.191,\ +3.319]$ & $[-0.167,\ +0.070]$ \\
threshold & $\frac{1}{45}$ & $0$ & $[-30.50,\ +0.005]$ & $[-7.417,\ +0.005]$ & $[-0.718,\ +0.005]$ \\
$\nu_{1/36}$ & $\frac{1}{36}$ & $0$ & $[-14.37,\ +0.212]$ & $[-3.248,\ +0.212]$ & $[-0.171,\ +0.212]$ \\
minimax, $\alpha \ge 1$ & $0.027537$ & $0$ & $[-15.03,\ +0.189]$ & $[-3.413,\ +0.189]$ & $[-0.189,\ +0.189]$ \\
minimax, $\alpha \ge \frac12$ & $0.031403$ & $0$ & $[-5.090,\ +1.048]$ & $[-1.048,\ +1.048]$ & $[-0.002,\ +0.710]$ \\
minimax, $\alpha \ge \frac13$ & $0.032615$ & $0$ & $[-2.438,\ +2.439]$ & $[-0.484,\ +1.947]$ & $[0,\ +0.885]$ \\
Paulson, $\nu_{\mathrm{WH}}$ & $\frac{1}{27}$ & $-\frac{1}{729}$ & $[-2.538,\ +3.866]$ & $[-0.124,\ +3.423]$ & $[0,\ +1.325]$ \\
\bottomrule
\end{tabular}
\end{table}

\paragraph{Comparison with Kerman's approximation.}
The improvement holds pointwise, not only in the worst case. On a grid covering
$\frac13 < \alpha \le \beta \le 10^6$, including points within $10^{-6}$ of the
diagonal in $\beta/\alpha$, the relative error of $\nu_{1/36}$ is everywhere smaller in magnitude
than that of $\nu_0$, and the ratio of the two stays below $\frac{13}{32} = 0.40625$.

\paragraph{Bounds.}
Although $\nu_{1/36}$ is an asymptotic upper bound on every ray (Section~\ref{sec:rays}), it is not an
upper bound everywhere: its negative errors are confined to a
region near the diagonal that shrinks as $\alpha$ grows.
The sign of the error
just off the diagonal changes at $\alpha \approx 2.059$, and, numerically on our grid at least, $\nu_{1/36} \ge \nuB$
throughout $2.06 \le \alpha \le \beta$.
Paulson's $\nu_{\mathrm{WH}}$ is also an asymptotic upper bound on
every ray, and on our grid it is an upper bound throughout $1 \le \alpha \le \beta$, but with six times
the worst error of $\nu_{1/36}$ there.

\paragraph{Smaller shape parameters.}
Below $\alpha = 1$ the errors of $\nu_{1/36}$ grow, mainly on the negative side (Table~\ref{tab:errors}).
Over $\alpha, \beta \ge \frac13$, its relative error is at most $14.4\%$ in magnitude, while Kerman's
approximation reaches $-100\%$ at $\alpha = \frac13$, where its $g$ vanishes; the minimax member for that
domain (Section~\ref{sec:otherc}) keeps the error within $2.44\%$.

\section{The ideal parameter}
\label{sec:idealc}

Setting $\nu_{c_1}(\alpha,\beta)$ equal to the median $m = \nuB(\alpha,\beta)$ and solving for $c_1$ gives
the value that makes the ratio form exact at $(\alpha,\beta)$,
\begin{equation}
c^*(\alpha,\beta) = \frac{m \left(\beta - \frac{1}{3}\right) - (1-m) \left(\alpha - \frac{1}{3}\right)}
{(1-m)/\alpha - m/\beta} \, .
\label{eq:idealc}
\end{equation}
It is symmetric: $c^*(\beta,\alpha) = c^*(\alpha,\beta)$. For $\alpha < \beta$, the denominator is positive,
since $m/(1-m) < 1 < \beta/\alpha$; on the diagonal, where every ratio form is exact, numerator and
denominator both vanish and $c^*$ is defined by continuity. By l'H\^opital's rule, with
$m' = \partial m/\partial\beta$ at $\beta = \alpha$,
\[
c^*(\alpha,\alpha) = \frac{\alpha^2\left(1 + 4m'\left(\alpha - \frac13\right)\right)}{1 - 4\alpha\, m'} \, .
\]
In general $m'$ has no closed form, so this does not give $c^*$ on the diagonal explicitly. At $\alpha = 1$,
however, $m = 1 - 2^{-1/\beta}$ gives $m' = -\frac12\log 2$, and so
\[
c^*(1,1) = \frac{3 - 4\log 2}{3\,(1 + 2\log 2)} \approx 0.03177 \, .
\]
Elsewhere on and near the diagonal, where \eqref{eq:idealc} loses accuracy to cancellation, the notebook
evaluates $c^*$ at a nearby pair symmetric about it. By the proof of Theorem~\ref{thm:cfamily},
$c^*(\alpha, r\alpha) \to c^*_r$ as $\alpha \to \infty$ along each ray.

Where $\alpha$ and $\beta$ both exceed $\frac16$, and $g_{c_1}$ is positive at both, $\nu_{c_1}(\alpha,\beta)$
increases with $c_1$ for $\alpha < \beta$: the derivative of $\log\big(g_{c_1}(\alpha)/g_{c_1}(\beta)\big)$
with respect to $c_1$ is $1/(\alpha g_{c_1}(\alpha)) - 1/(\beta g_{c_1}(\beta))$, and
$x\,g_{c_1}(x) = x^2 - x/3 + c_1$ increases for $x > \frac16$. In that region, then, the contour
$c^* = c_1$ is the curve on which $\nu_{c_1}$ is exact, and $\nu_{c_1}$ overestimates the median where
$c^* < c_1$ and underestimates it where $c^* > c_1$.

\begin{figure}[t]
  \centering
  \includegraphics{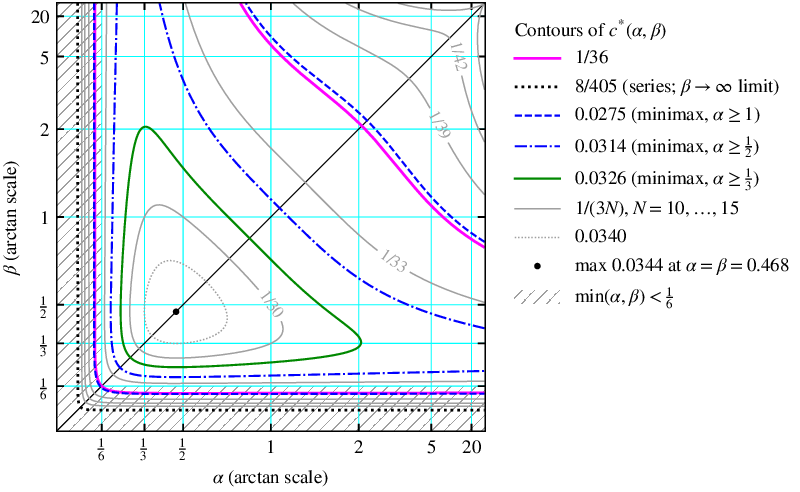}
  \caption{Contours of the ideal parameter $c^*(\alpha,\beta)$ of \eqref{eq:idealc}, on $\arctan$ axes
  covering $[0,\infty]^2$, with the diagonal drawn. Each contour is the curve on which the corresponding
  $\nu_{c_1}$ is exact; outside the shaded strips, $\nu_{c_1}$ overestimates the median on the side where
  $c^* < c_1$. Within the shaded strips, $\min(\alpha,\beta) < \frac16$, the contours for $c_1 < \frac{1}{36}$
  follow the smaller zero of $g_{c_1}$ and do not mark useful approximations. The dot marks the peak of
  $c^*$, and the dotted contour $c^* = 0.034$ shows its shape.}
  \label{fig:idealc}
\end{figure}

Figure~\ref{fig:idealc} shows contours of $c^*$ on $\arctan$ axes covering $[0,\infty]^2$. Numerically, $c^*$ has a single peak, on the diagonal: $c^* \approx 0.03440$ at $\alpha = \beta \approx 0.468$. Outside the shaded strips,
any $c_1$ above this maximum, such as $\frac{1}{29} \approx 0.03448$, therefore makes $\nu_{c_1}$ an upper
bound for all $\alpha < \beta$. Over
$\alpha, \beta \ge 1$, the surface lies between $0.0318$, at $\alpha = \beta = 1$, and $\frac{8}{405}$, approached as
$\alpha \to \infty$ with $\beta/\alpha \to \infty$; along the diagonal it decreases toward $\frac{1}{45}$.

\paragraph{The small-$c^*$ end.}
Toward the axes, $c^*$ falls to zero, but this part of the surface does not point to useful
approximations. For small $\alpha$, the median is extremely small, roughly like $2^{-1/\alpha}$, so an
exact ratio form needs $g_{c_1}(\alpha) \approx 0$.

\section{Other choices of \texorpdfstring{$c_1$}{c}}
\label{sec:otherc}

The minimax value of $c_1$ over $1 \le \alpha \le \beta$ is $0.027537 \approx 1/36.31$, with worst relative error
$0.189\%$ (Table~\ref{tab:errors}). It depends slightly on the error measure: minimizing the worst logit
error instead gives $0.027996 \approx 1/35.72$, and minimizing the worst deviation
$|I_{\nu_{c_1}}(\alpha,\beta) - \frac12|$ of the beta distribution function at the approximation gives
$0.027578 \approx 1/36.26$. So $\frac{1}{36}$ lies within the range of these minimax values. In terms of
Figure~\ref{fig:idealc}, each minimax value corresponds to the contour of $c^*$ that balances the worst errors
on its two sides.

Extending the domain below $\alpha = 1$ pushes the minimax value of $c_1$ up: over $\alpha \ge \frac12$ it
is $0.031403 \approx 1/31.84$, with worst error $1.05\%$, and over $\alpha \ge \frac13$ it is
$0.032615 \approx 1/30.66$, with worst
error $2.44\%$, attained at $\alpha = \frac13$, both near $\beta = 2\alpha$ and as $\beta \to \infty$.
Below $\alpha = \frac13$, no single $c_1$ does well.

A second coefficient can do somewhat better for larger $\alpha$: the series \eqref{eq:gammaseries}
truncated after its $c_2$ term has worst error $0.167\%$ over $\alpha \ge 1$, but $3.3\%$ over
$\alpha \ge \frac12$ and $30\%$ over $\alpha \ge \frac13$, and lacks the simple factored form.

In the upper right quadrant of Figure~\ref{fig:idealc}, where $\alpha, \beta \ge 1$, $c^*$ decreases as
$\alpha$ and $\beta$ grow, and most of that decrease goes with $1/\alpha + 1/\beta$, which suggests
letting $c_1$ depend on that sum. For example, $c_1 = \frac{8}{405} + \frac{184}{25515}\big(\frac{1}{\alpha} + \frac{1}{\beta}\big)$, which as
$\beta \to \infty$ reproduces the gamma median expansion \eqref{eq:gammaseries} through its $x^{-2}$ term,
has worst relative error $0.070\%$ over $1 \le \alpha \le \beta$, a third of that of $\nu_{1/36}$, and the
simpler $c_1 = (3 + 1/\alpha + 1/\beta)/150$ has $0.041\%$, a fifth. The result is no longer of the form
$g(\alpha)/(g(\alpha) + g(\beta))$, since $c_1$, and with it each of $g(\alpha)$ and $g(\beta)$, depends on
both $\alpha$ and $\beta$, and both choices are much worse below $\alpha = 1$, where $c_1$ grows like $1/\alpha$.

\section{Gamma median approximation}
\label{sec:gamma}

Each $g_{c_1}$ is itself an approximation to the gamma median: it agrees with the expansion
\eqref{eq:gammaseries} in its first two terms, $\alpha - \frac13$, and differs from it by $O(1/\alpha)$. Our recommended member is
\begin{equation}
g_{1/36}(\alpha) = \frac{(\alpha - \frac16)^2}{\alpha} = \alpha - \frac13 + \frac{1}{36\alpha} \, .
\label{eq:nuG3}
\end{equation}
By \eqref{eq:gammaseries},
$g_{1/36}(\alpha) - \nuG(\alpha) \sim \big(\frac{1}{36} - \frac{8}{405}\big)/\alpha = 13/(1620\,\alpha)$, so
it is an upper bound for large $\alpha$. Numerically it is an upper bound for all $\alpha > 0.7959$, where it
crosses $\nuG$, and for $\alpha \ge 1$ it overestimates by at most $0.212\%$, at $\alpha \approx 1.18$. At
$\alpha = 1$, it gives $\frac{25}{36}$ in place of $\log 2$, $0.19\%$ high.

\begin{figure}[t]
  \centering
  \includegraphics{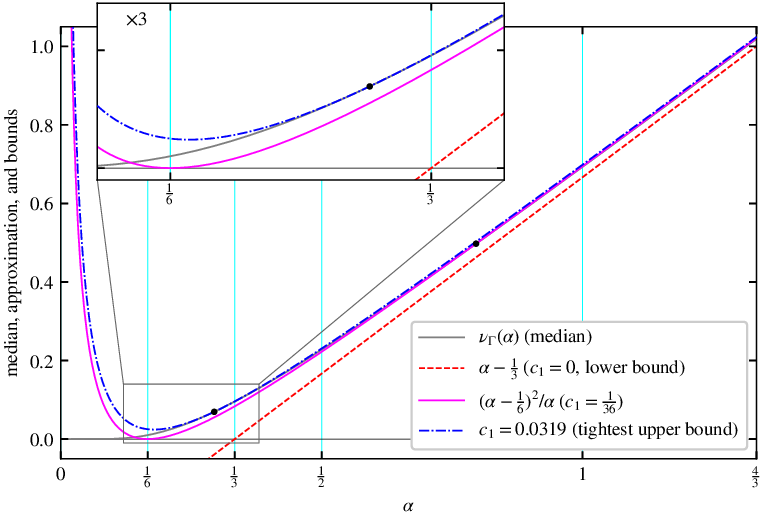}
  \caption{The gamma median $\nuG(\alpha)$ (gray) with the lower bound $\alpha - \frac13$ (red dashed), the
  approximation $(\alpha - \frac16)^2/\alpha$ (magenta), and the upper bound $\alpha - \frac13 + c_1/\alpha$ with
  $c_1 \approx 0.03194$, the smallest $c_1$ that makes it an upper bound for all $\alpha > 0$ (blue dash-dot). The inset magnifies the
  region $0.12 < \alpha < 0.38$ by a factor of 3. Dots mark where the upper bound is tangent to the median ($\alpha \approx 0.294$) and where the
  approximation crosses it to become an upper bound ($\alpha \approx 0.796$).}
  \label{fig:gammalin}
\end{figure}

Figure~\ref{fig:gammalin} compares it with $\nuG$, with the lower bound $\alpha - \frac13$, and with the
upper bound $\alpha - \frac13 + c_1/\alpha$ with $c_1 \approx 0.03194$ described below, for
$0 < \alpha \le \frac43$. Above about $\alpha = \frac12$, it is hard to distinguish from the median. Below, its
double zero at $\alpha = \frac16$ gives a soft parabolic floor, where the true median, roughly like
$2^{-1/\alpha}$, flattens toward zero. A larger $c_1$ lifts the floor: numerically, the smallest $c_1$ for which
$\alpha - \frac13 + c_1/\alpha \ge \nuG(\alpha)$ for all $\alpha > 0$ is $c_1 \approx 0.03194$, for which the two
touch near $\alpha \approx 0.294$. For small shape parameters, where the median behaves like $2^{-1/\alpha}$,
bounds and approximations built around that factor do much better \citep{lyon2021closed, lyon2023tight}.

\section{Numerical methods and use of AI tools}
\label{sec:numerics}

Medians were computed with the double-precision inverse incomplete beta and gamma functions of SciPy
\citep{virtanen2020scipy}. For each lower limit $\alpha_0$ of Table~\ref{tab:errors}, worst errors were taken
over $\alpha_0 \le \alpha \le 1000$ (400 logarithmically
spaced values) and $1 \le \beta/\alpha \le 10^5$ (400 values), with $\beta \le 10^6$, together with the
limit $\beta \to \infty$ computed from the gamma median. The minimax
values of $c_1$ were found by bounded scalar minimization of the worst error on this grid. The pointwise
comparison with Kerman's approximation used a grid denser near $\alpha = \frac13$ and near the diagonal:
$400$ linearly spaced values of $\alpha$ from $\frac13 + 10^{-4}$ to $0.6$ and $400$ logarithmically spaced
values from $0.6$ to $1000$; $60$ values of $\beta/\alpha$ with $\beta/\alpha - 1$ logarithmically spaced
from $10^{-6}$ to $10^{-1}$ and $300$ logarithmically spaced from $1.1$ to $10^6$; and $\beta \le 10^6$. Near the diagonal at large $\alpha$, where both errors approach rounding level, the
comparison was confirmed at selected points in 40-digit arithmetic with mpmath \citep{mpmath}, as was the
threshold in Remark~\ref{rem:threshold}. A Python notebook that reproduces the figures and the numbers quoted here is included
as an ancillary file with the arXiv version of this paper.

This note is a spin-off from a longer manuscript on bounds for the beta median, which the author wrote,
including most of its proofs, initially without AI tools; in later work on that manuscript, ChatGPT
(OpenAI) was useful in exposing errors in candidate proofs. In preparing the
note, the author used Claude (Anthropic) to help refactor and condense material from that manuscript, to
check algebra symbolically and numerically, and to port the author's Matlab scripts to the Python
notebook. Claude and Gemini (Google) were used to search for workable proof methods, to check and assess two
independently derived sets of proofs, and to synthesize the best of both into Appendix~\ref{app:proofs}.
Google's Deep Review tool \citep{google2026deepreview} reviewed several generations of the manuscript. The author checked the proofs, revised the text
throughout, and is responsible for its content. All numbers quoted can be reproduced from the notebook.

\section{Summary}

Replacing Kerman's $g(x) = x - \frac13$ by $g(x) = (x - \frac16)^2/x$, that is, adding $1/(36x)$, keeps
the ratio form simple, symmetric and exact on the diagonal, and reduces the worst relative error over
$1 \le \alpha \le \beta$ from $3.8\%$ to $0.21\%$. Numerically on our grid, and asymptotically by
Theorem~\ref{prop:ratio}, the new approximation is more accurate than Kerman's wherever the latter is
defined, with at most $\frac{13}{32}$ of its error. Theorem~\ref{thm:cfamily} places
$c_1 = \frac{1}{36}$ above the threshold $\frac{1}{45}$ for an asymptotic upper bound on every ray, the
ideal-parameter map of Section~\ref{sec:idealc} shows where it overestimates and underestimates, and
the same $g$ approximates the gamma median within $0.212\%$ for shape parameter at least $1$. For
posteriors under Kerman's neutral prior, the member $c_1 \approx 0.0326$ keeps the relative error within
$2.44\%$ over the extended domain $\frac13 \le \alpha \le \beta$.

\appendix
\section{Proofs}
\label{app:proofs}

\subsection{The median logit}

The proofs work with the logit, $\logit p = \log(p/(1-p))$, of the median. Three properties make it
convenient. Because the logit is increasing, it maps the median of a distribution on $(0,1)$ to the median
of the transformed distribution: $\logit\nuB(\alpha,\beta)$ is the median of $\logit X$ for
$X \sim \mathrm{Beta}(\alpha,\beta)$, whose density is smooth and log-concave on the whole real line. The
logit of a ratio form separates into a difference,
$\logit\big(g(\alpha)/(g(\alpha) + g(\beta))\big) = \log g(\alpha) - \log g(\beta)$, so it can be compared
with the median one shape parameter at a time. And the symmetry $\nuB(\beta,\alpha) = 1 - \nuB(\alpha,\beta)$
becomes a change of sign, with $\logit\nuB = 0$ on the diagonal.

Write $\Delta = 1/\alpha - 1/\beta$, a measure of the departure from the diagonal: it is zero there and
positive for $\alpha < \beta$. Along a ray $\beta = r\alpha$, we have $\Delta = (1 - 1/r)/\alpha$, so
approximations to the logit of the median can be ordered by powers of $1/\alpha$. We use three, of orders
$0$, $1$ and $2$:
\begin{align}
L_0(\alpha,\beta) &= \log\frac{\alpha}{\beta} \, , \qquad L_1(\alpha,\beta) = L_0(\alpha,\beta) - \frac{\Delta}{3} \, ,
\nonumber \\
L_2(\alpha,\beta) &= L_1(\alpha,\beta)
- \frac{\Delta\,(29\alpha^2 + 50\alpha\beta + 29\beta^2)}{810\,\alpha\beta\,(\alpha+\beta)} \, .
\label{eq:largeexp}
\end{align}
The zeroth, $L_0$, is the logit of the mean $\alpha/(\alpha+\beta)$. The first, $L_1$, agrees to first order with the logit of Kerman's
approximation, $\log(\alpha - \frac13) - \log(\beta - \frac13) = L_1 + O(\alpha^{-2})$, which is why that
approximation is correct to first order. The second, $L_2$, adds the second-order correction that decides the bounds of
Theorem~\ref{thm:cfamily}. Lemma~\ref{lem:largeexp} shows that it matches the logit of the median to second
order; its coefficients come out of Laplace's method in the proof.

\begin{lemma}[Large-shape expansion]
\label{lem:largeexp}
Fix $r \ge 1$. As $\alpha \to \infty$ with $\beta = r\alpha$,
\begin{equation}
\logit \nuB(\alpha,\beta) = L_2(\alpha,\beta) + O(\alpha^{-3}) \, .
\label{eq:logitexp}
\end{equation}
\end{lemma}

\begin{proof}
Let $n = \alpha + \beta$, $p = \alpha/n$
and $q = \beta/n$, which are fixed along the ray. The logit $W = \logit X$ of $X \sim \mathrm{Beta}(\alpha,\beta)$
has density proportional to $e^{\alpha w}(1+e^w)^{-n} = e^{n\varphi(w)}$, with $\varphi(w) = pw - \log(1+e^w)$,
and its median is $\logit\nuB$. With $\sigma(w) = 1/(1+e^{-w})$, we have $\varphi' = p - \sigma$ and
$\varphi'' = -\sigma(1-\sigma) < 0$, so $\varphi$ is strictly concave with its maximum at
$w^* = \log(p/q) = \log(\alpha/\beta)$, where $\sigma = p$. Put $\chi(x) = \varphi(w^*+x) - \varphi(w^*)$.
Differentiating $\sigma' = \sigma(1-\sigma)$ repeatedly gives
\[
\chi(x) = -\frac{A}{2}x^2 + \frac{B}{6}x^3 + \frac{C}{24}x^4 + \frac{D}{120}x^5 + O(x^6) \, ,
\]
with $A = pq$, $B = -pq(q-p)$, $C = -pq(1-6pq)$ and $D = -pq(q-p)(1-12pq)$.

The median $x_m$ of $W - w^*$ satisfies $\int_{-\infty}^{x_m} e^{n\chi} = \int_{x_m}^\infty e^{n\chi}$, that is,
\[
2\int_0^{x_m} e^{n\chi(x)}\,dx = \int_0^\infty \left(e^{n\chi(x)} - e^{n\chi(-x)}\right) dx \, .
\]
We expand both sides by Laplace's method \citep{olver1974asymptotics}. By concavity,
$\chi(x) \le \chi(\pm\delta)\,|x|/\delta$ for $|x| \ge \delta$, so the parts of the integrals with
$|x| \ge \delta$ are exponentially small in $n$. For $|x| < \delta$, put $x = y/\sqrt n$. Then
$n\chi = -Ay^2/2 + \varepsilon(y)$, where
$\varepsilon = By^3/(6\sqrt n) + Cy^4/(24n) + Dy^5/(120 n^{3/2}) + O(y^6 n^{-2})$, and for small $\delta$,
$|\varepsilon| \le Ay^2/4$, which gives the usual Gaussian domination of the remainders. Expanding
$e^\varepsilon$, the terms of order $n^{-1/2}$ and $n^{-3/2}$ are odd in $y$ and those of order $1$, $n^{-1}$
and $n^{-2}$ are even, so only the odd terms survive on the right. With
$\int_0^\infty y^{2j+1}e^{-Ay^2/2}\,dy = j!\,2^j/A^{j+1}$, the right side is $n^{-1/2}$ times
\[
\frac{2B}{3A^2}\,n^{-1/2} + \left(\frac{2D}{15A^3} + \frac{2BC}{3A^4} + \frac{16B^3}{27A^5}\right)n^{-3/2}
+ O(n^{-5/2}) \, ,
\]
where the $n^{-3/2}$ coefficient collects the terms $Dy^5/120$, $BCy^7/144$ and $B^3y^9/1296$ of
$e^\varepsilon$. On the left the integrand is close to $1$ near $0$, so $y_m = \sqrt n\,x_m = O(n^{-1/2})$;
then $\varepsilon = O(n^{-2})$ on $[0, y_m]$, and the left side is $n^{-1/2}$ times $2y_m - Ay_m^3/3 + O(n^{-5/2})$. Matching powers of $n$ gives
$x_m = \mu_1/n + \mu_3/n^2 + O(n^{-3})$, with
\[
\mu_1 = \frac{B}{3A^2} \, , \qquad
\mu_3 = \frac{D}{15A^3} + \frac{BC}{3A^4} + \frac{49B^3}{162A^5} \, ,
\]
where $\frac{49}{162} = \frac{8}{27} + \frac{1}{162}$ combines half the last coefficient on the right with
the term $A\mu_1^3/6$ from the left.

Since $(q-p)/(npq) = (\beta - \alpha)/(\alpha\beta) = \Delta$, the first term is
$\mu_1/n = -(q-p)/(3npq) = -\Delta/3$. Using $(q-p)^2 = 1 - 4pq$,
\[
\mu_3 = \frac{q-p}{p^2q^2}\left[-\frac{1-12pq}{15} + \frac{1-6pq}{3} - \frac{49(1-4pq)}{162}\right]
= -\frac{(q-p)(29 - 8pq)}{810\,p^2q^2} \, ,
\]
so $\mu_3/n^2 = -\Delta(29 - 8pq)/(810\,npq)$. As $n^2(29 - 8pq) = 29\alpha^2 + 50\alpha\beta + 29\beta^2$ and
$n^3pq = \alpha\beta(\alpha+\beta)$, this is the last term of \eqref{eq:largeexp}. Finally
$\logit\nuB = w^* + x_m$, and $O(n^{-3}) = O(\alpha^{-3})$ for fixed $r$.
\end{proof}

As a check, as $\beta \to \infty$ the expansion becomes $\log(\alpha/\beta) - 1/(3\alpha) - 29/(810\alpha^2)$,
which is $\log(\nuG(\alpha)/\beta)$ to second order by \eqref{eq:gammaseries}, consistent with
$\beta\,\nuB(\alpha,\beta) \to \nuG(\alpha)$. Numerically, $\alpha^3$ times the remainder settles to a constant
along each ray, about $0.0025$ for $r = 2$.

The expansion can also be obtained in the original variable. With $\psi(u) = -p\log(u/p) - q\log((1-u)/q)$,
the integrand $u^{\alpha-1}(1-u)^{\beta-1}$ is proportional to $e^{-n\psi(u)}/(u(1-u))$, and Temme's quadratic
change of variable $\frac12 z^2 = \psi(u)$, with $z$ of the sign of $u - p$, makes the exponent exactly Gaussian;
Watson's lemma then expands the two sides of the balance equation. This is the route of Temme's uniform
asymptotic expansions of the incomplete beta function \citep{temme1992asymptotic}. It gives the median
directly,
\[
\nuB = p - \frac{q-p}{3n} + \frac{q-p}{pq\,n^2}\Big(c^*_r - \frac{2pq}{9}\Big) + O(n^{-3}) \, ,
\]
which is the expansion of the ratio form $\nu_{c_1}$ with $c_1 = c^*_r$.

\subsection{Proof of Theorem~\ref{thm:cfamily}}

\thmbounds*

\begin{proof}
Write $a = 1/\alpha$ and $b = 1/\beta = a/r$, so $\Delta = a - b$. For large $x$,
$\log g(x) = \log x - \frac{1}{3x} + \big(c_1 - \frac{1}{18}\big)\frac{1}{x^2} + O(x^{-3})$, so along the ray
\[
\log\frac{g(\alpha)}{g(\beta)} = \log\frac{\alpha}{\beta} - \frac{\Delta}{3}
+ \Big(c_1 - \frac{1}{18}\Big)(a^2 - b^2) + O(\alpha^{-3}) \, .
\]
In terms of $a$ and $b$, the last term of \eqref{eq:largeexp} is
$-\Delta\,(29a^2 + 50ab + 29b^2)/(810(a+b))$. Subtracting from \eqref{eq:logitexp}, and using
$45(a+b)^2 - (29a^2 + 50ab + 29b^2) = 16(a+b)^2 + 8ab$,
\[
\logit\nuB - \log\frac{g(\alpha)}{g(\beta)}
= \Delta\left[\Big(\frac{8}{405} - c_1\Big)(a+b) + \frac{4ab}{405\,(a+b)}\right] + O(\alpha^{-3}) \, .
\]
With $b = a/r$ we have $ab/(a+b)^2 = r/(1+r)^2$, so the bracket is $(a+b)(c^*_r - c_1)$, and
\begin{equation}
\logit\nuB - \log\frac{g(\alpha)}{g(\beta)} = a^2(1 - r^{-2})(c^*_r - c_1) + O(a^3) \, .
\label{eq:logiterr}
\end{equation}
For large $\alpha$, both $g(\alpha)$ and $g(\beta)$ are positive, so the ratio form lies in $(0,1)$ and its
logit is $\log(g(\alpha)/g(\beta))$. It is below the logit of $\nuB$ when $c_1 < c^*_r$ and above it when
$c_1 > c^*_r$, for all sufficiently large $\alpha$, and the logit is increasing.
\end{proof}

\subsection{Proof of Theorem~\ref{prop:ratio}}

\thmerror*

\begin{proof}
By \eqref{eq:logiterr}, the logit error $\delta = \logit\nu_{c_1} - \logit\nuB$ is
$-a^2(1 - r^{-2})(c^*_r - c_1) + O(a^3)$. Since $\sigma'/\sigma = 1 - \sigma$, Taylor's theorem gives
$\sigma(w + \delta)/\sigma(w) - 1 = (1 - \sigma(w))\,\delta + O(\delta^2)$, so with $w = \logit\nuB$ the relative
error is $e_{c_1} = (1 - \nuB)\,\delta + O(\alpha^{-4})$. Along the ray, $\nuB = 1/(1+r) + O(a)$, so
$1 - \nuB = r/(1+r) + O(a)$, and $\frac{r}{1+r}(1 - r^{-2}) = 1 - \frac1r$, which gives the stated form. The
ratio of the values at $c_1 = \frac{1}{36}$ and $c_1 = 0$ therefore tends to $(c^*_r - \frac{1}{36})/c^*_r$.
As $r \to 1$, $c^*_r \to \frac{9}{405}$ and its magnitude tends to $\frac{2.25}{9} = \frac14$; as
$r \to \infty$, $c^*_r \to \frac{8}{405}$ and it tends to $\frac{3.25}{8} = \frac{13}{32}$; in between it is
monotone in $c^*_r$, which decreases as $r$ increases, since $r/(1+r)^2$ does.
\end{proof}

\section*{Acknowledgments}

The author thanks Jouni Kerman for helpful and actionable comments on a draft.
The author also acknowledges Gemini (Google), Claude (Anthropic), and ChatGPT (OpenAI)
for help with proof methods, symbolic and numerical checking, code translation,
refactoring of the draft, and critical reviews of proofs and drafts; their use is described in
Section~\ref{sec:numerics}.

\end{document}